
\def\date{\ifcase\month\or January\or February \or March\or April\or May %
\or June\or July\or August\or September\or October\or November           %
\or December\fi\space\number\day, \number\year}                          %
\magnification=\magstep1   
\hsize=15truecm            
\vsize=20truecm            
\parindent=0pt             
\line{\hfill\date}         
\vskip2\baselineskip       
\catcode`\@=11             
\newfam\msbfam             
\newfam\euffam             
\font\tenmsb=msbm10        
\font\sevenmsb=msbm7       
\font\teneuf=eufm10        
\font\seveneuf=eufm7       
\textfont\msbfam=\tenmsb   
\scriptfont\msbfam=\sevenmsb%
\textfont\euffam=\teneuf   
\scriptfont\euffam=
\seveneuf                  
\def\C{{\fam\msbfam C}}    
\def\N{{\fam\msbfam N}}    
\def\T{{\fam\msbfam T}}    
\def\SS{{\fam\msbfam S}}   
\def\R{{\fam\msbfam R}}    
\def\Z{{\fam\msbfam Z}}    
\font\bfone=cmbx10         
scaled\magstep1            
\catcode`\@=12             
\def\.{{\cdot}}            
\def\<{\langle}            
\def\>{\rangle}            
\def\({\big(}              
\def\){\big)}              
\def\hat{\widehat}         
\def\tor{\mathop{\rm tor}  
\nolimits}                 
\def\id{\mathop{\rm id}    
\nolimits}                 
\def\card{\mathop{\rm card}
\nolimits}                 
\def\pr{\mathop{\rm pr}    
\nolimits}                 
\def\defi{\buildrel\rm def 
\over=}                    
\def\ssk{\smallskip}       
\def\msk{\medskip}         
\def\bsk{\bigskip}         
\def\cen{\centerline}      
\def\implies
{\hbox{$\Rightarrow$}}     
\def\iff
{\hbox{$\Leftrightarrow$}} 
\font\smc=cmcsc10          
\def\subsetneq{\;\rlap{\raise3pt\hbox{${\subset}$}}%
{\raise-3pt\hbox{${\neq}$}}\;} 
\def\2{{\bf2}}             
\def\0{{\bf0}}             
\def\1{{\bf1}}             
\def\lead{\leaders\hbox to 1.5ex{\hss${.}$\hss}\hfill}
\def\arr{\hbox to 30pt{\rightarrowfill}}
\def\larr{\hbox to 30pt{\leftarrowfill}}

\def\mapright#1{\smash{\mathop{\arr}\limits^{#1}}}
\def\lmapright#1{\smash{\mathop{\arr}\limits_{#1}}}

\def\mapup#1{\Big\uparrow\rlap{$\vcenter{\hbox{$\scriptstyle#1$}}$}}
\def\lmapup#1{\llap{$\vcenter{\hbox{$\scriptstyle#1$}}$}\Big\uparrow}

\long\def\alert#1{\parindent2em\smallskip\line{\hskip\parindent\vrule%
\vbox{\advance\hsize-2\parindent\hrule\smallskip\parindent.4\parindent%
\narrower\noindent#1\smallskip\hrule}\vrule\hfill}\smallskip\parindent0pt}

\def\bx#1{\lower4pt\hbox{
\vbox{\offinterlineskip
\hrule
\hbox{\vrule\strut\hskip1ex\hfil{\uvertspace #1}\hfill\hskip1ex}
\hrule}\vrule}}
\newcount\litter                                         
\def\newlitter#1.{\advance\litter by 1                   
\edef#1{\number\litter}}                                 
\def\qedbox{\hbox{$\rlap{$\sqcap$}\sqcup$}}              
\def\qed{\nobreak\hfill\penalty250 
\nobreak\hfill\qedbox\vskip1\baselineskip\rm}            
\def\Cl{\hbox{$\cal C$}}
\def\lw{w_0}

\newlitter\comfort.
\newlitter\ferher.
\newlitter\flor.
\newlitter\glick.
\newlitter\hern.
\newlitter\hewross.
\newlitter\layer.
\newlitter\compbook.
\newlitter\leptin.
\newlitter\montzip.
\newlitter\reid.
\newlitter\ribes. 
\newlitter\scott.
\newlitter\walker.
\newlitter\zel.

\vglue 3\baselineskip
\cen{\bf Editorial Page}
\vskip2\baselineskip
\cen{\bfone The weights of closed subgroups of a locally compact group}
\msk
\cen{by Salvador Hern\'andez, Karl H. Hofmann, and Sidney A. Morris}
\bsk
{Salvador Hern\'andez}

{hernande@mat.uji.es}

{Universitat Jaume I, INIT and Depto de  Matem\'aticas,
Campus de Riu Sec, 12071 Castell\'on, Spain}
\msk
{Karl H. Hofmann, Corresponding author}

{hofmann@mathematik.tu-darmstadt.de}

{Fachbereich Mathematik, Technische Universit\"at Darmstadt, Schlossgartenstrasse 7, 64289 Darmstadt, Germany. }

 \msk
{Sidney A. Morris}

{morris.sidney@gmail.com}

{School of Science, IT, and Engineering, University of Ballarat,
Victoria 3353, Australia, and School of Engineering and Mathematical Sciences,
La Trobe University, Bundoora, Victoria 3086, Australia}

\vglue 2\baselineskip

{\bf Abstract.}\quad Let $G$ be an infinite locally compact group
and $\aleph$ a cardinal satisfying $\aleph_0\le\aleph\le w(G)$
for the weight $w(G)$ of $G$. It is shown that there is a closed
subgroup $N$ of $G$ with $w(N)=\aleph$. Sample consequences are:

{\parindent=2em
\item{(1)} Every infinite compact group contains an
infinite closed metric subgroup.
\item{(2)} For a locally compact group $G$ and $\aleph$
a cardinal satisfying $\aleph_0\le\aleph\le \lw(G)$, where $\lw(G)$ is the local
weight of $G$, there are either no infinite compact subgroups at all
or there is a compact  subgroup $N$ of $G$ with $w(N)=\aleph$.
\item{(3)} For an infinite abelian group $G$ there exists a properly
    ascending family
    of locally quasiconvex group topologies on  $G$, say,
    $(\tau_\aleph)_{\aleph_0\le \aleph\le \card(G)}$, such that
    $(G,\tau_\aleph)\hat{\phantom{m}}\cong\hat G$.

Items (2) and (3) are shown in Section 5.
}

\vfill\eject

\vglue2\baselineskip

\cen{\bfone The weights of closed subgroups of a locally compact group}
\msk
\cen{by Salvador Hern\'andez, Karl H. Hofmann, and Sidney A. Morris}
\bsk

{\bf Abstract.}\quad Let $G$ be an infinite locally compact group
and $\aleph$ a cardinal satisfying $\aleph_0\le\aleph\le w(G)$
for the weight $w(G)$ of $G$. It is shown that there is a closed
subgroup $N$ of $G$ with $w(N)=\aleph$. Sample consequences are:

{\parindent=2em
\item{(1)} Every infinite compact group contains an
infinite closed metric subgroup.
\item{(2)} For a locally compact group $G$ and $\aleph$
a cardinal satisfying $\aleph_0\le\aleph\le \lw(G)$, where $\lw(G)$ is the local
weight of $G$, there are either no infinite compact subgroups at all
or there is a compact  subgroup $N$ of $G$ with $w(N)=\aleph$.
\item{(3)} For an infinite abelian group $G$ there exists a properly
    ascending family
    of locally quasiconvex group topologies on  $G$, say,
    $(\tau_\aleph)_{\aleph_0\le \aleph\le \card(G)}$, such that
    $(G,\tau_\aleph)\hat{\phantom{m}}\cong\hat G$.

Items (2) and (3) are shown in Section 5.
}
\bigskip

\cen{\bf Introduction}

\msk
The {\it weight} $w(X)$ of a topological space $X$ is the
smallest cardinal $\aleph$ for which there is a basis
$\cal B$ of the topology of $X$ such that
$\card({\cal B})=\aleph$.
A compact group $G$ is {\it metric} iff its
weight $w(G)$ is countable, that is,
$w(G)\le\aleph_0$.
(See e.g. [\compbook], A4.10 ff., notably, A1.16.)
In particular, all compact Lie groups are metric.
It is not clear a priori that a compact group of uncountable weight contains
an infinite closed metric subgroup. Indeed, in Example 5.3(v) below we
will show that the precompact topological group defined on $\Z$
endowed with the Bohr-topology, which it inherits from its universal
almost periodic compactification, has no nonsingleton metric subgroup,
while  its weight is the cardinality of the continuum.

\msk
However, we shall prove the following  theorem which among other things will
show that every infinite compact group has an
infinite metric subgroup.

\msk

{\bf Main Theorem} \quad {\it Let $G$ be a locally compact group of uncountable
weight and let $\aleph_0\le\aleph<w(G)$. Then $G$ has a closed
subgroup $N$ with $w(N)=\aleph$.}
\msk
In other words, {\it for any infinite locally compact group $G$,
the entire interval
of cardinals $[\aleph_0,w(G)]$ is occupied by the weights of
closed subgroups of $G$.}

\ssk If $G$ is compact and connected, we shall see that
$[\aleph_0,w(G)]$ is filled
even with  closed {\it  normal} and indeed connected subgroups.
It remains unsettled
whether a profinite group has, in this sense, enough {\it normal}
 closed subgroups.

\msk
We shall deal with a proof of the Main Theorem  in a piecemeal way.

For reasons of presenting a stepwise proof
let us call \Cl\ the class of all Hausdorff topological groups $G$
satisfying the following condition
\msk
{\it for each infinite cardinal $\aleph\le w(G)$ there is a closed
subgroup $H$ of weight $\aleph$.}
\msk
Our Main Theorem says that all locally compact groups are
contained in \Cl.

We first aim to show that all compact groups are in \Cl\ and we begin
with compact abelian groups.

\bsk
\cen{\bf1. Compact abelian groups}
\msk

For an abelian group $A$, let $\tor A$ denote its torsion subgroup.

The first portion of our first observation is a consequence of a more precise
statement due to {\smc W. R. Scott} [\scott]. Since we prove what we shall
need in a shorter way (and quite differently) we present a proof which
will also establish the second part. A divisible hull of an abelian
group we are using can
be constructed in a very special way by using the results of
{\smc E. A. Walker} in [\walker].

\msk
\bf
Proposition 1.1. \quad {\it Let $A$ be an
uncountable abelian group and $\aleph_0\le\aleph<|A|$. Then
$A$ contains a subgroup $B$ such that $(A:B)=\aleph$.

Moreover,
If $(A:\tor A)$ is at least $\aleph$
then $B$ may be picked so that $B$
is pure and $A/B$ is torsion free.}
\msk
\bf
Proof.\quad\rm Let $D$ be a divisible hull of $A$ according to
[\compbook], Proposition A1.33.
 If $D=D_1\oplus D_2$  is any direct decomposition, and
 $\pr_1\colon D\to D_1$ is the projection
onto the first summand of $D$,
then $A/(A\cap D_2)\cong (A+D_2)/D_2\cong \pr_1(A)$.
Now any subgroup of $D_1$ is a subgroup of the divisible hull of $A$
and therefore meets $A$ and thus $A\cap D_1\subseteq\pr_1(A)$ nontrivially
(see [\compbook], Proposition A1.33); therefore $D_1$ is a divisible
hull of $\pr_1(A)$. Hence either $\pr_1(A)$ is finite or else
$\card D_1=\card\pr_1(A)=\card(A/(A\cap D_2))$ by [\compbook],
Proposition A1.33(i).
\msk
Since we control $\card D_1$ by choosing $D_1$ appropriately,
we aim to set $B=A\cap D_2$ and thereby prove our first assertion.
We thus have to exclude the possibility that $pr_1(A)$
might turn out to be finite by an inappropriate choice of $D_1$.
We now let $\tor A$ denote the torsion subgroup  of $A$.
Then $\tor D$ is a divisible hull of $\tor A$, and
$D\cong(\tor D)\times (D/\tor D)$ by [\compbook], p.~ 657,
Proposition A1.38.
 We now  distinguish two cases:
\ssk
(a) Case $\card(\tor A)=\card A$.
Since $A$ is uncountable, one of the $p$-primary components of
$\tor A$, as $p$
ranges through the countable set of primes,
say $A(p)$, satisfies $\card A(p)=\card(\tor A)=\card A$.
In particular $\card A(p)$ is uncountable, that is,
its $p$-rank $\card A$ is uncountable and
agrees with the $p$-rank of $D(p)$ (see [\compbook], p.~656,
Corollary A1.36(iii). In view of $D(p)\cong(\Z(p^\infty)^{(\card A)})$
by [\compbook], p.~659, Theorem A1.42(iii), we find a direct
summand $D_1$ of $D(p)$ of $p$-rank $\aleph$, giving us a
direct summand of $\tor D$ and thus yielding a
direct sum decomposition $D=D_1\oplus D_2$. Since the $p$-rank
$\aleph$ of $D_1$ is infinite, and $D_1$ is the divisible
hull of $\pr_1(A)$ we know that $\pr_1(A)$ cannot be finite,
whence $\aleph=\card D_1=\card(A/(A\cap D_2)$. Our first
assertion then follows with $B=A\cap D_2$.
\ssk
(b) Case $\card(A/\tor A)=\card A$. Then the (torsion free) rank
of $D$ is $\card A$ (see [\compbook], p.~656, Corollary A1.36(iii)),
 By the structure theorem of
divisible groups (see e.g.\ [\compbook], Theorem A1.42]) and
elementary cardinal arithmetic, we can write
$D=D_1\oplus D_2$  with a torsionfree subgroup $D_1$ of
cardinality $\aleph$. Then $\pr_1(A)\subseteq D_1$ cannot be
finite, and as in the first case, we let
 $B=A\cap D_2$ and have $\aleph=\card D_1=\card A/B$
as in our first assertion.
\bsk
It remains to inspect the case that
$\card (A/\tor A)\ge\aleph$. Then we may assume
$\tor A\subseteq D_2$ and $D_1$ torsion free; but then
$\tor A\subseteq A\cap D_2=B$ whence $A/B$ is torsion free.\qed

The second part of the preceding proposition is also a consequence
of Theorem 4 of [\walker].

\msk

\bf Corollary 1.2. \quad \it Every uncountable abelian group has a
proper subgroup
of index $\aleph_0$. \qed

\msk
(A different, but likewise not entirely trivial proof by Hewitt and Ross
is found in [\hewross], p.\ 227.)

As usual, for a topological group $G$, the identity component of $G$
will be denoted $G_0$.
\msk\bf
Corollary 1.3.\quad \it Let $G$ be an
infinite  compact abel\-ian group and assume $\aleph_0\le\aleph<w(G)$. Then
$G$ has a closed subgroup $M$ with $w(M)=\aleph$.

Moreover, if $w(G_0)\ge \aleph$, then $M$ may be chosen to
be connected. \qed
\msk
\bf
Proof. \quad \rm  By the Annihilator Mechanism
([\compbook], Theorem 7.64), for a subgroup $H$ of $G$
and its annihilator $H^\perp$ in $\hat G$,
one has $\hat H\cong\hat G/H^\perp$ and thus
$w(H)=\card \hat H = \card  \hat G/H^\perp$.
Since $\card \hat G=w(G)$ we have $\card \hat G>\aleph$.

Moreover, $H$ is connected iff $\hat G/H^\perp$ is
torsion free. (See [\compbook], Corollary 8.5.)
The Corollary  is therefore equivalent to  the
Proposition 1.1 above. \qed
\msk
\bf
Corollary 1.4. \quad \it Let $\aleph$ be an
infinite cardinal and let  $G$ be a compact
connected group with
$w(G)>\aleph$.   Then $G$ has a compact connected abelian
subgroup $T$ with $w(T)=\aleph$.
\msk
\bf
Proof.\quad\rm Let $T$ be a maximal compact connected abelian
subgroup of $G$.  Then $w(T)=w(G)$ by [\compbook], Theorem 9.36(vi).
Then the assertion follows from Corollary 1.3. \qed

Let us see what we have
in the case of $\aleph=\aleph_0$.
A  compact abelian group $A$ is metric iff $\hat A$ is
countable ([\compbook], Theorem 7.76). Thus
Corollary 1.3 trivially implies
\msk\bf
Corollary 1.5. \quad {\it Every compact abelian group has
an infinite closed metric subgroup.}
\bsk\rm
A result by {\smc Efim Zelmanov} says:
\msk
\bf
Theorem 1.6.\quad \it An infinite compact group contains an infinite abelian
subgroup.
\msk
\bf
Proof.\quad \rm  See [\zel] or [\ribes], p. 162.\qed
\msk
This together with Corollary 1.5, implies
\msk

\bf
Corollary 1.7.\quad \it  Every compact group contains
an infinite metric compact abelian subgroup. \qed

\rm

\msk While this corollary answers the question whether infinite compact groups
have infinite compact metric subgroup in the affirmative,
 we should keep in mind, that {\smc Zelmanov}'s Theorem in
itself is not a simple matter. It therefore appears worthwhile to
pursue the question further.
\bsk

\cen{\bf 2. Connectivity versus total disconnectivity in compact groups}
\msk
\bf
Lemma 2.1.\quad\it Let $G$ be an arbitrary compact group.
Then the following conclusions hold:

{\parindent2em

\item{\rm(i)} $G$ and $G_0\times G/G_0$ are homeomorphic.
In particular, if $G$ is infinite, then
$$w(G)=\max\{w(G_0),w(G/G_0)\}.$$

\item{\rm(ii)} There is a profinite subgroup $D$ of $G$
such that $G=G_0D$ and $G_0\cap D$ is normal in $G$ and central
in $G_0$.  If $w(G_0)< w(G)$, then $w(G)=w(G/G_0)=w(D)$.

\item{\rm(iii)} If $G$ is profinite, then $G$ and
$(\Z/2\Z)^{w(G)}$ are homeomorphic.

}

\msk
\bf Proof. \quad \rm For (i) see [\compbook], 10.38.

Regarding (ii), see [\compbook], 9.41 and note
$G/G_0\cong D/(D\cap G_0)$.

For (iii), see [\compbook], 10.40. \qed

After Lemma 2.1 the question whether a compact group $G$ is in $\Cl$
 splits into two cases:

Case 1. If $w(G_0)=w(G)$, one may assume that $G$ is connected.

Case 2. If $w(G_0)<w(G)$, one may assume that $G$ is profinite.

\msk
We recall from Corollary 1.4 that for a compact {\it connected group}
$G$ the set of all infinite cardinals $\le w(G)$ is filled with  the
set of all infinite cardinals representing the weights of closed
connected abelian subgroups. In the following we amplify this observation

\msk\bf
Proposition 2.2. \quad \it Let $\aleph$ be an infinite
 cardinal such that $\aleph<w(G)$
for a  compact connected group $G$. Then $G$ contains a closed
{\rm connected and normal}
subgroup $N$ such that $w(N)=\aleph$.
\msk
\bf
Proof.
\quad \rm Following the Levi-Mal'cev Structure Theorem for Compact Connected
Groups ([\compbook], Theorem 9.24) we have $G=G'Z_0(G)$ where
the algebraic commutator subgroup $G'$ is a characteristic
compact connected semisimple subgroup and the identity component of the
center $Z_0(G)$ is a characteristic compact connected abelian subgroup.
\msk
Case 1. $w(G')\le \aleph$. Then $w(G)=w(Z_0(G))$, and by
Corollary 1.3, $Z_0(G)$ contains a connected closed subgroup $N$
of weight $\aleph$; since it is central, it is normal.
\msk
Case 2. $w(G')>\aleph$. If we find a compact connected normal subgroup
$N$ of $G'$, we are done, since the normalizer of $N$ contains both
$G'$ and the central subgroup $Z_0(G)$, hence all of $G=G'Z_0(G)$.
Thus it is no loss of generality to assume that $G=G'$ is a compact
connected semisimple group. \ssk

Case 2a. $G=\prod_{j\in J}G_j$  for a family of compact connected
(simple) Lie groups. Then $w(G_j)=\aleph_0$, and $\aleph< w(G)
=\max\{\aleph_0,\card J\}$ (see e.g.\ [\compbook], EA4.3.). Since
$\aleph$ is infinite and smaller than $w(G)$, we have $w(G)=\card J$.
Then we find a subset $I\subseteq J$ such that $\card I=\aleph$,
and set $N=\prod_{i\in I}G_i$. Then $w(N)=\card N=\aleph$.
\ssk
Case 2b. By the Sandwich Theorem for Semisimple Compact
Connected Groups ([\compbook], 9.20)  there is a family of
simply connected compact simple Lie groups $S_j$ with
center $Z(S_j)$ and there are surjective morphisms
$$\prod_{j\in J}S_j\mapright{f}G\mapright{q}\prod_{j\in J}S_j/Z(S_j)$$
such that $qf$ is the product $\prod_{j\in J}p_j$ of the
quotient morphisms $p_j\colon S_j\to S_j/Z(S_j)$.
Now both products $\prod_{j\in J}S_j$ and $\prod_{j\in J}S_j/Z(S_j)$
have the same weight $\card J$ which agrees with the
weight of the sandwiched group $G$. Define $I$ as in Case 2a
and set $N=f(\prod_{i\in I}S_i)$ and note that
$q(N)=\prod_{i\in I}S_i/Z(S_i)$. Hence $N$ is sandwiched between
two products with weight $\card I=\aleph$ and hence has weight
$\aleph$. This proves the existence of the asserted $N$ in the
last case. \qed

\bsk

\cen{\bf 3. The generating degree}
\msk
Now we have to reach beyond connectivity, all the while still staying within
the class of compact groups.

We refer to a cardinal invariant for compact groups $G$ which is one
of several alternatives to the weight $w(G)$, namely, the so called
{\it generating degree} $s(G)$ (see [\compbook], Definition 12.15).
The definition
relies on the Suitable Set Theorem, loc. cit. Theorem 12.11, which
in turn invokes the so called Countable Layer Theorem
(see [\layer] or [\compbook], Theorem 9.91).
Indeed recall that in a compact group $G$ a subset $S$ is called {\it suitable}
iff it does not contain $1$, is closed and discrete in $G\setminus\{1\}$,
and satisfies $G=\overline{\<S\>}$.
The Suitable Set Theorem asserts, that every compact group $G$
has a suitable set.
A suitable set is called  {\it special} iff its cardinality
is minimal among all suitable subsets of $G$. The cardinality $s(G)$
of one, hence every special suitable set
is called {\it the generating degree of} $G$.

\msk
The relevance of the generating degree in our context is the following

\msk
\bf Proposition 3.1. \quad\it Let $G$ be a profinite group with uncountable
weight. Then $w(G)=s(G)$.

\msk \bf
Proof.\quad \rm By Proposition 12.28 of [\compbook], for an infinite
profinite, that is, compact totally disconnected group we have
$$w(G)=\max\{\aleph_0,s(G)\}.$$
This implies the assertion immediately in the case of $w(G)>\aleph_0$.
\qed

The next step, namely, proving that every profinite
group is in $\Cl$
will be facilitated by a lemma on suitable sets for which
all ingredients are contained in [\compbook].

\msk
\bf Lemma 3.2. \quad \rm (a) \it Let $S$ be any suitable set of a compact group
$G$. Then $\card S\le w(G)$.

{\rm(b)} If $G$ is profinite and $S$ is an infinite suitable subset of $G$,
then $\card S=w(G)$.

\ssk

\bf Proof.\quad\rm
(a) Let $\cal B$ be a basis of the topology of $H$
of cardinality $w(H)$. Since $S$ is discrete in $G\setminus\{1\}$,
for every element $x\in S$ there
is an element $U(x)\in{\cal B}$ with $U(x)\cap X=\{x\}$.
Then $x\mapsto U(x):S \to {\cal B}$ is an injective function
and thus $\card S\le \card{\cal B}=w(H)$.
(Cf. [\compbook], p.~620, proof of 12.16.)
\ssk
(b) Assume  that $G$ is profinite and  $w(G)$ is uncountable. Then
$w(G)=s(G)\le \card S$ by Proposition 3.1  and the definition
of $s(G)$. From this and (a), $w(G)=\card S$ follows.
Now assume that   $S$ is infinite. Then either
$S$ is uncountable or $\card S=\aleph_0$. In the first case,
 $w(G)$ is  uncountable by (a) and $w(G)=\card S$ holds.
 If, on the other hand,
$\card S=\aleph_0$ and $w(X)=\aleph_0$, then $\card S=w(X)$ as well.  \qed

The significance of this lemma is that

\ssk
{\it for a profinite group,
infinite suitable subsets
{\rm all}  have the same cardinality, namely, the weight of
the group}.
\msk

It is instructive  to take note of the following remarks
which are pertinent to this context:
\msk

\bf Remark a. \quad\rm The universal monothetic and the universal solenoidal
compact groups $G$ have generating degree $s(G)=1$, density $d(G)=\aleph_0$
(i.e. they are separable), and  weight $w(G)=2^{\aleph_0}$. (For the
concept of {\it density} see e.g.\ [\compbook], p.~620, Definition 12.15.)
\msk

\bf Remark b.\quad \rm If
 $H$ is a precompact group whose Weil completion $G$ is profinite and
if $H$ has an infinite relatively compact suitable subset,
then $\card H\geq w(H)$.

\msk
\bf
Proof. \quad \rm Assume that $S$ is an infinite relatively compact
 suitable subset of $H$.
Then by [\compbook], p.~616, Lemma 12.4,
$S$ is a suitable subset of $G$ and  by Lemma 3.2 it follows that
$\card H\geq \card S= w(G)=w(H)$. \qed

\msk
\bf Remark c. \quad \rm
(i)  For every compact group $G$ and infinite cardinal number $\aleph$,
the inequalities $\aleph<w(G)\leq 2^\aleph$ imply $d(G)<w(G)$.
\ssk

(ii) Let $G$ be a compact group  of weight $\aleph_1$.
Then $G$  contains  countable dense
subgroups. If $G$ is profinite, then none of these  contains an
infinite relatively compact suitable set.

\msk
{\bf Proof.} \quad (i) This follows from the equation
$d(G)=\log w(G)$ valid for any compact group $G$,
see [\comfort].

(ii) Let $G$ be a compact group of weight $\aleph_1$. Then
by (i) it has a countable dense subgroup $H$. Suppose
that $G$ is profinite and $H$ has an infinite relatively compact
suitable
subset $S$. Then by
Remark b we would have $\aleph_0=\card H\ge w(H)=w(G)=\aleph_1$,
a contradiction. \qed

\bsk

Now we show that every profinite group of uncountable weight is in $\Cl$.
\msk
\bf Lemma 3.3.\quad \it
Let $G$ be a profinite group of uncountable weight and let $\aleph<w(G)$
be an infinite cardinal. Then there is a closed subgroup $H$
such that $w(H)=\aleph$.

\msk
\bf
Proof. \quad \rm Let $T$ be a suitable subset of $G$ with
$\card T=w(G)$  according
to Proposition 3.1. Then $T$ contains a subset $S$ of cardinality
$\aleph$. We set $H=\overline{\<S\>}$. Now $S$ is discrete in
$H\setminus\{1\}$ since $T$ is discrete in $G\setminus\{1\}$.
Hence $S$ is an infinite suitable subset of the profinite
group $H$. Hence, by Lemma 3.1(b), $w(H)=\card (S)=\aleph$ follows.\qed

\msk
\bf
Corollary 3.4. \quad\it An infinite profinite group contains an infinite
compact metric subgroup.
\msk
\bf
Proof.\quad  \rm
Let $G$ be an infinite metric group.
If $w(G) =\aleph_0$ then $G$ itself is metric. If $G$ has uncountable weight,
then we apply Lemma 3.3 with $\aleph=\aleph_0$. \qed

\bsk
Now we are ready to prove that every compact group is in $\Cl$ which
is the main portion of the following result:

\msk
\bf
Theorem 3.5. \quad\it Let $\aleph$ be an infinite cardinal and $G$
a compact group such that $\aleph< w(G)$. Then there is a closed subgroup
$H$ such that $w(H)=\aleph$, and if $G$ is connected, $H$ may be chosen
normal and connected.
\msk
                                                                  \bf
Proof.\quad \rm Let $G_0$ denote the identity component of $G$.
The case that $w(G)=w(G_0)$ is handled in Proposition 2.2. So we assume
$w(G_0)<w(G)$. By Lemma 2.1 we may assume that $G$ is totally disconnected,
that is, profinite. Then Lemma 3.3 proves the assertion of the theorem. \qed
\msk
In particular, we have the following conclusion:
\msk\bf
Corollary 3.6. \quad \it Every infinite compact group contains an infinite
 closed metric subgroup.   \qed\rm

Recall that by Corollary 1.7 we know that we find even a closed {\it abelian}
metric subgroup. For a compact connected group $G$, Corollary 1.4 shows
that for any infinite cardinal $\aleph\le w(G)$ there is in fact a closed
connected abelian subgroup of weight $\aleph$.

\msk
\bf
Problem.\quad\it If $G$ is a compact group and $\aleph$ is an infinite
 cardinal $\le w(G)$. Is there a {\rm normal} closed subgroup $H$
such that $w(H)=\aleph$? \rm
\msk
The answer is affirmative if it is affirmative for profinite groups.
A profinite group has weight $\aleph$ iff its set of open-closed
normal subgroups has cardinality $\aleph$. This observation
points into the direction of an affirmative answer to the problem.

\bsk

\cen{\bf 4. The weights of closed subgroups of a locally compact group}
\bsk

Now we finish the proof of the Main Theorem by showing that every
{\it locally} compact group is in $\Cl$. A first step is the
following:
\msk
{\bf Lemma 4.1.} \quad\it Every  locally compact pro-Lie group belongs
to \Cl. \rm
\msk
{\bf Proof.}\quad Let $N$ be a compact normal subgroup
of $G$ such that $G/N$ is a Lie group. Then $w(G/N)=\aleph_0$.
We conclude $w(G)=w(N)$. We apply Theorem A to $N$ and obtain
a closed subgroup of $N$ of weight $\aleph$.\qed
\msk
In Montgomery's and Zippin's classic [\montzip] we note Lemma 2.3.1 on
p.~54 and the Theorem on p.~175 and thus find
\msk

{\bf Lemma 4.2.} \quad\it Every locally compact group contains an open
pro-Lie group. \qed

\msk
{\bf Lemma 4.3.} \quad \it Let $H$ be an infinite
open subgroup of a topological
Hausdorff group $G$. Let $\cal X$ be a set of cosets $Hg$ of $G$
modulo $H$. Then $w(\<\bigcup{\cal X}\>)=\max\{w(H),\card({\cal X}))$.

\msk
\bf Proof.\quad\rm Abbreviate $\<\bigcup{\cal X}\>$ by $K$.
Then $w(H)\le w(K)$ and $\card({\cal X})\le w(K)$ and thus
$$\max\{\card({\cal X}), w(H)\}\le w(K).\leqno(1)$$

In order to prove the reverse inequality, let
 $D$ be a dense subset of $H$ of cardinality $d(H)$. For
each finite tuple ${\cal F}=(Hg_1,Hg_2,\dots,Hg_n)\in {\cal X}^n$,
$n\in\N$, the set $P_{\cal F}\defi Hg_1Hg_2\cdots Hg_n$ has a dense
subset $\Delta_{\cal F}\defi Dg_1Dg_2\cdots Dg_n$ of cardinality
$d(H)$ with the density of $H$.
Then $K=\bigcup_{\cal F}P_{\cal F}$, with the union extended
over the set of finite $n$-tuples $\cal F$,  has a dense subset
$\Delta=\bigcup_{\cal F}\Delta_{\cal F}$ whose cardinality is
$\le \card{\cal X}\.d(H)\le \max\{\card({\cal X}), w(H)\}$
(see [\compbook], p.~620, Proposition 6.20 and its proof).
It follows that
$$w(K)\le \card(\Delta)\.w(H)\le \max\{\card({\cal X}), w(H)\}.
\leqno(2)$$
Now (1) and (2) together prove the assertion.\qed
\msk
The following conclusion now completes the proof of the
Main Theorem:
\msk
{\bf Corollary 4.4.}\quad \it Every locally compact group
belongs to \Cl.

\msk
{\bf Proof.}\quad \rm If $\aleph=w(G)=\aleph_0$ there is nothing to prove.
Assume $\aleph_0\le \aleph<w(G)$; we have to find a closed
subgroup such that $w(K)=\aleph$. By Lemma 4.2, let $H$
be an open pro-Lie subgroup of $G$.
Now $w(G)=\max\{w(H),\card(G/H)\}$. If $w(G)=w(H)$, we find
$K$ by Lemma 4.1. If $w(G)=\card(G/H)$ then we pick a subset
$\cal X$ of $G/H$ of cardinality $\aleph$. Then
$K=\<\bigcup{\cal X}\>$ has weight $\card{\cal X}=\aleph$
by Lemma 4.3. \qed

\bsk

\cen{\bf 5. An application of the Main Theorem}
\bsk
\noindent Let $(G,\tau )$ be an arbitrary topological abelian group.
We shall continue to write abelian groups additively, unless specified
otherwise such as in the case of the multiplicative circle group
$\SS^1=\{z\in\C: |z|=1\}$.
A character on $(G,\tau )$ is a continuous morphism
$\chi\colon G\to \T=\R/\Z$. The pointwise sum of two
characters is again a character, and the set $\widehat{G}$ of all characters
is a group with pointwise multiplication as the composition law. If
$\widehat{G}$ is equipped with the compact open topology $\widehat \tau$,
it becomes a
topological group $(\widehat{G},\widehat{\tau })$ which is called the
{\it dual group} of $(G,\tau )$.

\ssk
{\bf
Definition 5.1} (Varopoulos) \quad Let $G$ and $H$ be
two abelian groups then
we say that {\it $G$ and $H$  are in duality}
if and only if there is a $\Z$-bilinear function

$$\<\.,\.\>\colon G\times H\longrightarrow \T=\R/\Z $$

such that
$$\eqalign{%
&(\forall 0_G\ne g\in G)(\exists h\in H)\, \<g,h\>\ne0_\T,\hbox{and}\cr
&(\forall 0_H\ne h\in H)(\exists g\in G)\, \<g,h\>\ne0_\T.\cr}$$

\msk

{\bf
Definition 5.2} (Varopoulos) \quad Assume that $G$ and $H$ are in
duality. Then {\it a topology $\tau $ on $G$ is compatible
with the duality} if $(G,\tau )\,\hat{}=H$.

\msk
Typically, if $G$ is a locally compact abelian group and $\hat G$ its
Pontryagin dual, that is, its character group, then $G$ and $\hat G$
are in duality, but there are group topologies $\tau$ on $G$ which are
in general coarser than the given locally compact topology on $G$ but
which are nevertheless compatible with this duality. One of the best known
is the so-called  {\it Bohr topology} $\tau^+$ on $G$ which we
discuss in the next example for
the sake of completeness and because it plays a role in our subsequent
discussion which, at least in the case that $G$ is discrete, provides
a substantial cardinality of topologies on $G$ which are compatible
with the Pontryagin duality of $G$.

\msk
{\bf
Example 5.3.} \quad Let $G=(A,\tau)$
be a locally compact abelian group and
let $\eta\colon G\to G^\alpha$ be the Bohr compactification morphism
and let $\tau^+$ be the pull-back topology on $A$ (that is, the topology
which makes $\eta$ an embedding $\epsilon\colon G^+\to G^\alpha$).
Equivalently,
$\tau^+$ is the topology of pointwise convergence when $A$ is considered
as the space of characters of $\hat G$ via Pontryagin duality.
We write $G^+$
for the topological group $(A,\tau^+)$. Then

{\parindent2em

\item{(i)} If $(a_n)_{n\in\N}$ is a sequence of $A$, and $a\in A$,
then $a=\tau\hbox{-}\lim_na_n$ iff $a=\tau^+\hbox{-}\lim_na_n$. If $\tau$ is
the discrete topology, then $(a_n)_{n\in\N}$ converges w.r.t.\
$\tau^+$ iff it is eventually constant.
\ssk
\item{(ii)} $\tau^+$ is compatible with the duality between $A$ and $\hat G$,
that is, $\hat{G^+}=\hat G$.  Also, $\hat{G^+}=(G^\alpha)\hat{\phantom{x}}$.
\ssk
\item{(iii)} If $\sigma$ is any group topology on $A$ such that
$\tau^+\subseteq\sigma\subseteq\tau$, then $\sigma$ is
compatible with the duality between $A$ and $\hat G$.
\ssk
\item{(iv)} $w(G)=w(\hat G)$, $w(G^+)=w(G^\alpha)=\card(\hat G)$.
\ssk
\item{(v)} If $G$ has no nonsingleton compact subgroups, such as
$G=\Z$ or $G=\R$, then $G^+$ has no nonsingleton metric subgroups.

}
\msk

\bf Proof. \quad\rm For easy reference we provide proofs.

(i) \quad The fact that $G$ and $G^+$ have the same converging sequences
is based on two implications of which
``$a=\tau\hbox{-}\lim_na_n$ implies $a=\tau^+\hbox{-}\lim_n a_n$'' is immediate
from the continuity of $\eta$ and the definition of $\tau^+$.
The other implication has fascinated several authors (see e.g.\
[\flor], [\glick], [\leptin], [\reid]);
the following argument was credited by {\smc Reid}
 to originate from Varopoulos in the 1960s but may have been
around before:
Let $a=\tau^+\hbox{-}\lim_n a_n$. We consider $G$ as the character
group of $\hat G$. For any $f\in L^1(\hat G)$ we calculate the value
of the Fourier transform $\hat f\in C_0(G)$ at $a_n$ as
$$\hat f(a_n)= \int_{\hat G}f(\chi)\exp(-i\<\chi,a_n\>)\,d\chi$$
By the definition of $\tau^+$ as topology of pointwise convergence,
the sequence of functions $f(\bullet)\exp(-i\<\bullet,a_n\>)$ is
dominated by $|f|$ since the absolute value of the exponential is 1,
and it converges pointwise to $f(\bullet)\exp(-i\<\bullet,a\>$.
Hence by the Lebesgue Dominated Convergence Theorem we have
$$\lim_n \hat f(a_n)=\hat f(a).$$
Since the algebra of Fourier transforms $\hat f$ is uniformly dense
in $C_0(G)$, this implies $a=\tau\hbox{-}\lim_n a_n$.

\ssk

(ii)\quad Let $\chi\colon G^+\to\T$ be a (continuous) character.
Then there is a unique character
$\overline{\chi}\colon G^\alpha\to \T$
such that $\overline{\chi}\circ\eta=\chi$ by the definition
of $\tau^+$. However, if $\phi\colon G\to\T$ is a character, then
there is a unique character $\phi^\alpha\colon G^\alpha\to \T$
such that $\phi^\alpha\circ\eta=\phi$ by the universal
property of the Bohr compactification. Hence there is a
bijective correspondence
$\beta\colon\hat G\to(G^\alpha)\hat{\phantom{x}}$. Then
$\chi\mapsto \beta^{-1}(\overline{\chi}):\hat{G^+}\to \hat G$
is a bijection such that
$\beta^{-1}(\overline{\chi})\circ\eta=\chi$.
The following commutative diagram may help
$$\matrix{G^+&\mapright{\epsilon}&G^\alpha\cr
\lmapup{\id_A}&&\mapup{\id_{G^\alpha}}\cr
G&\lmapright{\eta}&G^\alpha.\cr}$$
Thus a character $A\to\T$ is $\tau^+$--continuous iff it is
$\tau$-continuous.

\ssk

(iii) \quad Every character $\chi\in\hat G$, $\chi\colon G=(A,\tau)\to\T$
is $\tau^+$--continuous by (ii). Then it is $\sigma$-continuous
since $\tau^+\subseteq\sigma$. Conversely, let
$\chi\colon (A,\sigma)\to\T$ be a $\sigma$-continuous character,
then it is $\tau$-continuous, since $\sigma\subseteq\tau$.
Hence a character $A\to\T$ is $\sigma$--continuous iff it is $\tau$
continuous.

\ssk

(iv) \quad For the equalities $w(G)=w(\hat G)$ and
$w(G^\alpha)= \card(G^\alpha)\hat{\phantom{x}}$
see e.g. [\compbook], Theorem 7.76(i) and (ii),
p.~364. Since $G^+$ has a dense homeomorphic image in $G^\alpha$
we have $w(G^+)=w(G^\alpha)$. Now  $(G^\alpha)\hat{\phantom{x}}
=\hat G$ by (ii) above. This proves the assertion.

\ssk

(v) \quad
 Let $H^+$ be a first
countable nonsingleton  subgroup of $G^+$.
 Then $H^+$ is precompact
since $G^+$ is precompact, and its topology is determined by its
convergent sequences since it is first countable. On the other hand,
 by (i) above, it has the same convergent sequences as
the subgroup $H$ which has the same underlying group as $H^+$
but whose topology is the one induced by the topology of $G$.
 Hence $H^+=H$ as topological groups. Since $H^+$ is
precompact, the locally compact group $\overline H$ is precompact
and thus is compact. Now assume that $G$ is a locally compact
abelian group without compact nondegenerate subgroups, say,
$G=\Z$ or $G=\R$. In those cases $G^+$ cannot have any
nonsingleton metric subgroup.  \qed

\msk By its very definition, the Bohr topology $\tau^+$ is a
precompact topology. Indeed,
 a topological group $G$ is said to be {\it precompact} or
{\it totally bounded} if    for
any neighborhood $U$ of the neutral element in $G$, there is
a subset $F\subseteq G$ with $\card(F)< \aleph_0$ such that
$FU=G$. Next  we need to generalize the concept of total boundedness:

\msk

{\bf Definition 5.4} \quad Let $\aleph$ be a cardinal number.
A topological group $G$ is said to be $\aleph$-{\it bounded}
when for every neighborhood $U$
of the neutral element in $G$, there is a subset $S\subseteq G$ with
$\card(S) < \aleph$ such that
$SU=G$. \qed

According to this definition, a group $G$ is totally bounded iff
it is $\aleph_0$-bounded. It is uniformly Lindel\"of if it is
$\aleph_1$-bounded.

For a topological abelian  group $\Gamma$ and a cardinal $\aleph$
we let ${\cal K}_\aleph(\Gamma)$ denote the set of all compact subsets
$K\subseteq\Gamma$ with $w(K)<\aleph$.

\msk
\bf
Lemma 5.5. \quad \it Let $G$ be the character group of an
abelian topological group $\Gamma$. Let $\aleph$ be a cardinal
$\le w(\Gamma)$, and let $G$ have the topology of uniform
convergence on compact subsets $K\in{\cal K}_\aleph(\Gamma)$.
Then $G$ is $\aleph$-bounded.

\msk

\bf
Proof. \rm \quad From Theorem 3.4 in [\ferher] we have the following
information due to Ferrer and  Hern\'andez who deal with the following
set-up:

Let $X$ be a set, let $M$ be a metrizable space,
and let $Y$ be a subset of $M^X$
that is equipped with some bornology $\cal B$
consisting of pointwise relatively compact sets.
Denote by $\mu_{\cal B}$  the uniformity on $X$ defined
as $\sup\{\mu_F\,:\,F\in{\cal B}\}$.
The concept of an $\aleph$-bounded topological group generalizes
rather immediately to that of an $\aleph$-bounded uniform space.
Now we have

\ssk
\noindent
{\bf The $\aleph$--Boundedness Theorem.}\quad
{\it If $\aleph$ is a cardinal such that
$w(M)<\aleph$, then
the following statements are equivalent:

{\parindent2em

\item{\rm(i)} $(\forall F\in{\cal B})\,w(F)<\aleph$.

\item{\rm(ii)} $(X,\mu_{\cal B})$ is $\aleph$-bounded.

}}

\ssk
Now
we apply this result with $\Gamma=X$, $\T=M$, $G=Y$,
${\cal K}_\aleph(\Gamma)={\cal B}$, and, finally, with
the topology of uniform
convergence on sets of ${\cal K}_\aleph(\Gamma)$ being the
uniform topology of $\mu_{\cal B}$. We have $w(M)=w(\T)=\aleph_0
\le \aleph$. Hence the implication (i) $\Rightarrow$ (ii)
of the $\aleph$--Boundedness Theorem yields
the assertion of the Lemma. \qed

\bsk

Let $G$ be a topological group, the {\it local weight} (or {\it character})
$\lw(G)$ is the smallest among the cardinals of  neighbourhood bases
at the neutral
element.  We say that a collection $\{K_i : i\in I\}$
is a compact cover of $G$ when $\cup_{i\in I} K_i=G$ and
$K_i$ is compact for all $i\in I$. The {\it compact covering number} $\kappa(G)$
of $G$ is defined as the smallest of the cardinals of the members of
the set of compact covers of $G$.

If $G$ is any locally compact group, let $H$ be an almost connected open
subgroup and let $C$ be a maximal compact subgroup of $H$.
Then $G$ is homeomorphic to $\R^n\times C\times G/H$ and
$\lw(G)=\max\{\aleph_0,w(C)\}$ and $\kappa(G)=\card G/H$. If $G$ is abelian and
$\Gamma=\hat G$, then $\lw(\Gamma) = \kappa(G)$ and
$\kappa(\Gamma)= \lw(G)$. In this sense, $\lw$ and $\kappa$ are ``dual''
cardinals.

In view of our Main Theorem we may summarize:

\msk
{\bf The Local Weight Lemma for Locally Compact Groups.}\quad \it
 For a locally
compact nondiscrete group $G$ select any almost connected open subgroup $H$ and
any maximal compact compact subgroup $C$ of $H$. Then

{\parindent2em
\item{\rm(1)} $\lw(G)=\max\{\aleph_0, w(C)\}$,

\item{\rm(2)} $w(G)=\max\{\lw(G), \card(G/H)\}$.

\item{\rm(3)} If $\aleph_0\le\aleph\le\lw(G)$, then either $G$ contains
no infinite compact subgroups or else there is a compact subgroup of
weight $\aleph$. \qed

}
\msk\rm

Consider now  a locally compact abelian group $G$ with dual group
$\Gamma$. If $\lw(\Gamma)$ is uncountable, then $\Gamma$ contains
a compact subgroup $C_\Gamma$ with weight $w(C_\Gamma)=\lw(\Gamma)$.
Now let $\aleph\leq \kappa(G)=\lw(\Gamma)$ be any infinite cardinal. Consider
$G$ as the character group of the locally compact group  $\Gamma$.

Write
 $\tau_\aleph$ for the the
topology of uniform convergence on the sets $K\in{\cal K}_\aleph(\Gamma)$.
We note that by definition,
\ssk
\noindent
{\it all topological abelian groups $(G, \tau_{\aleph_{\alpha}})$,
$0\le\alpha\in\Omega$ are locally quasiconvex.}
\ssk
Let $\Omega$ the initial set of ordinals $0<1<\cdots\le\xi$ such that
$\{\aleph_\alpha:\alpha\in\Omega\}$ is the interval $[\aleph_0,\lw(\Gamma)]$
of cardinals.1
\ssk
It is in the proof of the following lemma that we use the principal result
of this paper.
\msk

\bf Main Lemma 5.6.\quad \it For $\alpha<\xi$ in $\Omega$, the following
conclusions hold
{\parindent2em
\item{\rm(i)} ${\cal K}_{\aleph_{\alpha}}\subsetneq
                    {\cal K}_{\aleph_{\alpha+1}}$,
\item{\rm(ii)} $\tau_{\aleph_\alpha}\subsetneq
                \tau_{\aleph_{\alpha+1}}$.

}

\msk
\bf
Proof. \rm (i)\quad By the Main Theorem (indeed by Theorem 3.5),
the compact group $C_\Gamma$  contains a
subgroup $K$ of  of weight $w(K)=\aleph_\alpha$. What is relevant
for the proof is the fact that $\Gamma$ contains a compact subspace $K$ of
weight $\aleph_\alpha$.
 Thus
$K\in{\cal K}_{\aleph_{\alpha+1}}\setminus{\cal K}_{\alpha}$.
The inclusion ${\cal K}_{\aleph_\alpha}(\Gamma)
\subseteq {\cal K}_{\aleph_{\alpha+1}}(\Gamma)$
is trivial.

(ii) \quad From (i) it follows that
$\tau_{\aleph_\alpha}\subseteq\tau_{\aleph_{\alpha+1}}$.
{}From  Lemma 5.5 we know that $(G, \tau_{\aleph_\alpha})$ is
$\aleph_\alpha$-bounded and
that $(G, \tau_{\aleph_{\alpha+1}})$ is $\aleph_{\alpha+1}$-bounded.
Suppose, by way of contradiction that the topologies
$\tau_{\aleph_\alpha}$ and $\tau_{\aleph_{\alpha+1}}$ were equal.
Then the implication (ii) $\Rightarrow$ (i) of the
$\aleph$--Boundedness Theorem would imply that
every compact subset $K\subseteq \Gamma$ of weight
$w(K)<\aleph_{\alpha+1}$ would have weight $w(K)<\aleph_{\alpha}$.
This would contradict (i) above, and this contradiction shows
that $\tau_{\aleph_{\alpha+1}}$ is strictly finer than $\tau_{\aleph_\alpha}$.
\qed

The Main Lemma 5.6 establishes the essential result of this section:

\msk
{\bf
Theorem 5.7.}
{\it Let $(G,\tau)$ be a locally compact abelian group which has uncountable compact covering number
$\kappa(G)$. Then

\ssk{\parindent2em

\item{\rm(i)} for each cardinal
$\aleph$ with $\aleph_0\le\aleph\le\kappa(G)$ there is a
locally quasiconvex group topology $\tau_\aleph$ on $G$ such that
$(G,\tau_\aleph)\hat{\phantom{x}}=\hat G$, and

\item{\rm(ii)} for $\aleph_0\le\aleph<\aleph'\le\kappa(G)$ one
has $\tau_\aleph\subsetneq\tau_{\aleph'}\subseteq \tau$. \qed

}}

Thus, if $[\aleph_0, \kappa(G)]$ denotes the full interval of
infinite cardinals up to the cardinality of $G$, then Theorem 4.7
provides $\card[\aleph_0,\kappa(G)]$-many
locally quasiconvex group topologies
on the abelian group $G$ all of which are coarser than the original topology
of $G$ and yield the (locally compact abelian ) group
$\hat G$  as character group. This means that all topologies
$\tau_\aleph$ have the same compact subsets as $\tau$.
\bsk
\msk

{\bf
Corollary 5.8.}
{\it Let $G$ be an uncountable abelian group. Then

\ssk

{\parindent2em

\item{\rm(i)} for each cardinal
$\aleph$ with $\aleph_0\le\aleph\le\card(G)$ there is a
locally quasiconvex group topology $\tau_\aleph$ on $G$ such that
$(G,\tau_\aleph)\hat{\phantom{x}}=\hat G$, and

\item{\rm(ii)} for $\aleph_0\le\aleph<\aleph'\le\card(G)$ one
has $\tau_\aleph\subsetneq\tau_{\aleph'}$. \qed

}}

Again, if $[\aleph_0, \card(G)]$ denotes the full interval of
infinite cardinals up to the cardinality of $G$, then Corollary 5.8
provides $\card[\aleph_0,\card(G)]$-many
locally quasiconvex group topologies $\tau$
on the abelian group $G$ all of which yield  the
(compact abelian) group
$\hat G$  as character group of $(G,\tau)$.
\bsk

\cen{\bf References}

\msk
{\parindent 1.5em

\item{[\comfort]} Comfort, W. W., {\it Topological Groups},
in: K. Kunen and J. E. Vaughan, Eds., North-Holland, Amsterdam, 1984,
 Chapter 24, 1143--1263.

\item{[\ferher]} Ferrer, Maria V., Hern\'andez, Salvador,
{\it Dual topologies on
non-abelian groups}, Submitted to the Proceedings of  the conference
``Algebra meets Topology,'' Barcelona, July 2010.

\item{[\flor]} Flor, P.,
{\it Zur Bohr-Konvergenz der Folgen,}
Math. Scand. {\bf 23} (1968), 169--170.

\item{[\glick]} Glicksberg, I.,
{\it Uniform boundedness for groups},
Canadian J. Math. {\bf 14} (1962), 269--276.

\item{[\hern]}
Hern\'andez, Salvador, {\it Questions raised at the
conference} \rm ``Algebra meets Topology,'' Barcelona, July 2010.
{\tt hernande@mat.uji.es}

\item{[\hewross]}
Hewitt, E., and K. A. Ross,
``Abstract Harmonic Analysis I,''
Grundlehren {\bf115} Springer Verlag, Berlin etc., 1963.

\item{[\layer]}  Hofmann, K. H., and S. A. Morris,
{\it A structure theorem on compact groups}, Math. Proc.
Camb. Phil. Soc. {\bf130} (2001), 409--426.

\item{[\compbook]} ---,
``The Structure of Compact Groups'', de Gruyter, Berlin, 2$^{\rm nd}$
Edition 2006,  xvii+858 pp.

\item{[\leptin]} Leptin, H.,
{\it  Abelsche Gruppen mit
kompakten Charaktergruppen und Dua\-li\-t\"ats\-theo\-rie
gewisser linear topologischer abelscher Gruppen},
Abh. Math. Sem. Univ. Hamburg {\bf 19} (1955), 244-263.

\item{[\montzip]} Montgomery, D., and Zippin, L,
``Topological transformation groups,''
Robert E. Krieger Publishing Co., 1955. xi+282 pp.
2nd edition 1966, xi+289 pp.

\item{[\reid]} Reid, G. A.,
{\it On sequential convergence in groups,}
Math. Z. {\bf 102} (1967), 225--235.

\item{[\ribes]} Ribes, L., and P. Zaleskii,
``Profinite Groups,''
Springer-Verlag, Berlin etc., 2000, xiv+435 pp.
2nd edition 2010, xvi+464 pp.

\item{[\scott]} Scott, W. R.,
{\it The number of subgroups of given index in a nondenumerable
group},
Proc. Amer. Math. Soc. {\bf5} (1954), 19--22.

\item{[\walker]} Walker, E. A.,
{\it Subdirect sums and infinite abelian groups},
Pac. J. Math. {\bf9} (1959), 287--291.

\item{[\zel]} Zelmanov, E. I.,
{\it On periodic compact groups}, Israel J. of Math. {\bf77} (1992), 83--95.

}

\bsk

{Salvador Hern\'andez}

{hernande@mat.uji.es}

{Universitat Jaume I, INIT and Depto de  Matem\'aticas,
Campus de Riu Sec, 12071 Castell\'on, Spain}
\msk
{Karl H. Hofmann}

{hofmann@mathematik.tu-darmstadt.de}

{Fachbereich Mathematik, Technische Universit\"at Darmstadt, Schlossgartenstrasse 7, 64289 Darmstadt, Germany. }
 \msk
{Sidney A. Morris}

{morris.sidney@gmail.com}

{School of Science, IT, and Engineering, University of Ballarat,
Victoria 3353, Australia, and School of Engineering and Mathematical Sciences,
La Trobe University, Bundoora, Victoria 3086, Australia}

\bye